\input amstex
\documentstyle{amsppt}
\NoBlackBoxes
\pageheight{530pt}
\pagewidth{330pt}
\def\bk{\overline{k}}
\def\bQ{\Bbb Q}
\def\bZ{\Bbb Z}

\def\proj{\operatorname{Proj}}
\def\tF{\tilde{F}}
\def\tE{\tilde{E}}

\def\cP{\Cal P}
\def\cX{\Cal X}
\def\cQ{\Cal Q}
\def\oR{\overline{R}}

\def\bP{\Bbb P}
\def\tE{\tilde{E}}
\def\af #1.{\Bbb A^{#1}}
\def\au#1.{\operatorname {Aut}\,(#1)}
\def\tC{\tilde{C}}
\def\tD{\tilde{D}}

\def\ring#1.{\Cal O_{#1}}
\def\pr #1.{\Bbb P^{#1}}
\def\pic#1.{\operatorname {Pic}\,(#1)}
\def\bbbq{{\Bbb Q}}

\def\dim{\operatorname{dim}}
\def\qpic#1.{{\pic #1.}_{\bbbq}}

\def\ses#1.#2.#3.{0\longrightarrow #1 \longrightarrow #2 \longrightarrow #3
\longrightarrow 0}
\def\Aut{\operatorname{Aut}}
\def\mg{\overline {M}_g}
\def\ag{\operatorname{A}_g}

\def\fmgn#1.{\overline{\Cal M}_{g,#1}}

\def\vmgn#1.#2.{\overline{M}_{#1,#2}}
\def\vmg #1.{\overline{M}_{#1}}
\def\mgn{\overline{M}_{g,n}}

\def\vmon#1.{\vmgn 0.#1.}
\def\ms#1.{\overline{\Cal M}_{#1}}
\def\vln#1.{\overline{L}_{#1}}

\def\fgn{\overline{F}_{g,n}}
\def\bmgn#1.{\partial{\mgn #1.}}
\def\ugn#1.{{\Cal U}_{g,#1}}
\def\sm#1.#2.{\overline{\Cal M}_{#1,#2}}
\def\su#1.#2.{{\Cal U}_{#1,#2}}

\def\ex#1.{{\Bbb E}(#1)}
\def\mc#1.{\overline{NE}_1(#1)}
\def\md#1.{\overline{NE}^1(#1)}
\def\tmc#1.{\overline{NE}^{T}_1(#1)}
\def\tpic#1.{\operatorname{Pic}^{T}(#1)}
\def\tqpic#1.{{\tpic #1.}_{\bbbq}}
\def\omc#1.{{NE}_1(#1)}
\def\tomc#1.{{NE}^{T}_1(#1)}
\def\rmc#1.#2.{\overline{NE}_1(#1/#2)}

\def\sym#1.#2.{\operatorname {Sym}^{#1}(#2)}
\def\rpic#1.#2.#3.{\operatorname {Pic}_{#1/#2}^{#3}}

\def\nsd#1.{\operatorname{NS}^{1}(#1)}

\define\til #1.{\tilde{#1}}
\def\uring#1.{\ring #1.^*}
\def\tX{\tilde{X}}
\def\tcX{\tilde{\cX}}

\def\bF{\Bbb F}

\nologo
\topmatter
\title Polarized pushouts over finite fields \endtitle
\author Sean Keel \endauthor
\address Department of Mathematics, The University of Texas, Austin, TX 78712
\endaddress
\email keel\@math.utexas.edu \endemail
\dedicatory To Steve Kleiman on his sixtieth birthday \enddedicatory
\thanks The author received partial 
support from Harvard University, a  University of
Texas Faculty Research Grant,  and NSF grant DMS-9988874. \endthanks
\abstract
Let $p:Y \rightarrow X$ 
be a surjection between schemes projective over the algebraic 
closure of a finite field. Let $L$ be a 
line bundle on $X$ such that $p^*(L)$ is globally generated. 
A natural necessary and sufficient condition is given under
which some positive tensor power of $L$ is globally generated. 
An application is a sufficient condition
for semi-ampleness of nef line bundles on 
${\overline{M}_{g,n}}$ in positive characteristic, which
implies the semi-ampleness of $\lambda$, $\psi_i$
and myriad other nef line bundles. 
\endabstract
\endtopmatter


\subhead\S 0 Introduction and statement of results \endsubhead

Here I consider the following question: Let
$p: Y \rightarrow X$ be a proper surjection. Suppose
$p^*(L)$ is globally generated. Under what 
conditions does it follow that $L$ is semi-ample
(ie. some positive tensor power is globally generated)?
This always holds when $X$ is normal, so the interesting
case is when $p$ is the normalisation. 

\definition{0.0 Definition} Let $L$ be a nef line 
bundle on a complete algebraic space $X$. 
Two closed points $x_1,x_2 \in X$
are {\bf $L$-equivalent} iff there is a connected closed
curve $x_1,x_2 \in C \subset X$ with $L \cdot C =0$.
We say that $L$-equivalence is {\bf bounded} if there is
an integer $m > 0$ so that any two $L$-equivalent points
lie on a connected curve $C$ as above, with $C$ having at most
$m$ irreducible components.
\enddefinition

The main technical result of this paper is the following:

\proclaim{0.1 Theorem} 
Let $p:\tX \rightarrow X$ be a surjection between
algebraic spaces proper over the algebraic closure
of a finite field. 
Let $L$ be a 
line bundle on $X$. $L$ is semi-ample iff $p^*(L)$
is semi-ample and $L$-equivalence is bounded.
\endproclaim

(0.1) fails for all other fields, take for 
example the normalisation of a rational curve
with one ordinary node, and $L$ a degree zero
line bundle. The
pullback is always trivial, but $L$ is semi-ample iff
torsion. One cannot drop the boundedness assumption
in (0.1), see \S 5. 

Here are some applications of (0.1).

\cite{Keel99} gives a general basepoint
freeness 
result in positive characteristic (that fails in
characteristic zero)
namely 
a nef line bundle $L$ is semi-ample iff its restriction
to the exceptional locus, $\ex L.$, semi-ample, where
$\ex L.$ is the union of $L$-exceptional 
subvarieties, and 
a subvariety $Z \subset X$ is called $L$-exceptional
if $c_1(L)^{\dim(Z)} \cdot Z =0$ (or equivalently,
$L|_Z$ is not big). 
The main difficulty
in applying this in practice is that $\ex L.$
can have bad singularities, e.g. it is frequently
reducible, so one would like to pass to a desingularisation.
An immediate corollary of (0.1) is:

\proclaim{0.3 Corollary} Let $L$ be a nef line bundle
on a scheme $X$ projective over the algebraic closure
of a finite field. Then $L$ is semi-ample iff the
restriction of $L$ to the normalisation of the
exceptional locus is semi-ample and $L$-equivalence
is bounded.
\endproclaim

(0.3) has
applications to $\mgn$.

\proclaim{0.3.1 Corollary} Characteristic $p > 2$. Let
$L$ be a nef line bundle on $\mgn$. $L$ is 
semi-ample iff $L$-equivalence is bounded, and
$r^*(L)$ is semi-ample, for
$$
r: \vmgn 0. 2g +n. \rightarrow \mgn
$$
the natural finite map whose image is the locus
of curves with only rational components.
\endproclaim

The boundedness condition is a kind of combinatorial
question. I expect it always holds, see (0.6), and
can show it does  in many cases. A precise statement requires
some notation:

For a subset $S \subset N :=\{1,2,3,\dots,n\}$
$\pi_S: \vmgn g. N. \rightarrow \vmgn g. S.$
is the tautological map dropping the points labeled by
$S^c$ (and stabilizing). Let $C$ be a stable pointed
curve. Varying the moduli of the normalisation
of one of the irreducible components (together
with its distinguished points), keeping the other
components fixed, induces a natural finite map: 
$h: \vmgn r. m. \rightarrow \mgn$ (where 
$r$ is the geometric genus of the varying component).
The product of these maps over all components of
$C$ has image the closed stratum of topological
type $C$. For details see \cite{GKM00}.

\proclaim{0.4 Corollary} Characteristic $p > 2$. Let
$L$ be a nef line bundle on $\mgn$. Assume for each
map $h: \vmgn r. m. \rightarrow \mgn$ as above with
$r \leq 1$ 
that either
$h^*(L)$ is numerically trivial, or 
$L = \pi_S^*(L_S)$ for a (necessarily nef)
line bundle $L_S$ on $\vmgn r. S.$ with
$\ex L_S. \subset \partial \vmgn r. S.$.
Then $L$ is semi-ample.
\endproclaim

(0.3.1) and (0.4) are proved in \S 3, see
(3.6). (0.4) covers myriad cases. Here are just
a few:

\proclaim{0.5 Corollary} In characteristic $p > 0$ the
line bundles $\lambda$, $\psi_i$, and $D_{g,n}$ of
\cite{GKM00,4.8}, are all semi-ample. \endproclaim

Recall $\lambda$ is the determinant of the Hodge
bundle (which has fibre at a curve the vector space
of global $1$-forms).
Semi-ampleness of $\lambda$ is known in all characteristics.
In characteristic $p > 0$ it is due to Szpiro (by
a completely different argument), \cite{Szpiro90}. 
The associated map
is to the Satake compactification of $\ag$. 

$\psi_i$ is the bundle with fibre the
cotangent bundle of the curve at the $i^{th}$ point.
It is nef and big in all characteristics, but
fails to be semi-ample in characteristic zero, see
\cite{Keel99}. 
I claimed the semi-ampleness for $p > 0$ in
\cite{Keel99}, but the proof contains a gap -- I
use implicitly the {\it fact} that for a product
of $\vmgn g_i. n_i.$ the boundary is connected, which
fails for $\vmgn 0. 4.$ (in fact this is the
only case where it fails).

$D_{g,n}$ is introduced in \cite{GKM00}.
The semi-ampleness of $D_{g,n}$ gives an interesting
birational contraction: 
Let $\fgn \subset \mgn$ be the locus of
flag curves, i.e. the image 
map $f: \vmgn 0. g + n./S_g \rightarrow \mgn$
induced by gluing on $g$ copies of a one pointed
irreducible rational curve with a single node 
(any two are isomorphic) at $g$ points (which the symmetric
group $S_g \subset S_{g+n}$ permutes). $f$ is the normalization
of $\fgn$. The flag locus is important from several points
of view, see e.g. \cite{HarrisMorrison98}, \cite{Logan00}.
It is particularly important in the
Mori theory of $\mgn$: By \cite{GKM00,0.7}, $D_{g,n}$ is nef  and
the corresponding face of the (Kleiman) Mori cone
is 
$$
D_{g,n}^{\perp} \cap \omc{\mgn}. = \omc {\vmgn 0. g+n.}/S_g. .
$$
The face is natural both geometrically and purely 
combinatorially, and outside of this face the Mori
cone is completely known.
The semi-ampleness of $D_{g,n}$ means this face 
can be blowndown.
\proclaim{0.5.1 Corollary} In positive characteristic 
there is a birational map 
$$
q: \mgn \rightarrow Q_{g,n}
$$ 
such that $q(C)$ is a point
iff $C$ is numerically equivalent to a curve in $\fgn$.
The relative Mori cone $\rmc \mgn. Q_{g,n}.$ is naturally identified
with $\mc {\vmgn 0. g+n.}/S_g.$. \endproclaim

By \cite{GKM00} the analog of the assumption
in (0.4) holds for any nef line bundle $L$ on
$\vmgn r. m.$ for $r \geq 2$. A slight generalization
holds for $r =1$, namely any $L$ is a tensor
product of $L_S$, and I believe (0.4) could be
generalized to allow for this. 
I do not know of a nef line bundle 
on $\vmgn 0. M.$ which does not satisfy this
generalized condition. 

I propose the following:

\proclaim{0.6 Conjecture} Every nef line bundle on $\mgn$
is semi-ample in positive characteristic. \endproclaim

Note that together with \cite{GKM00,0.2} this gives
the conjecture:

\proclaim{0.7 Conjecture} The Mori cone of 
$\mgn$ is generated by finitely many rational curves
and in characteristic
$p>0$ every face has an associated contraction. 
\endproclaim

Thus in characteristic $p>0$ the Mori cone behaves
(conjecturally) like that of a Fano variety, although 
$\mgn$ is usually of general type.

\proclaim{0.8.1 Question} If $L$ is a non-trivial nef
line bundle on $\vmgn 0. N.$ can it be written as
a tensor product of line bundles
of the form $\pi_S^*(L_S)$ for $L_S$ as in
(0.4)? 
\endproclaim

It is natural to wonder quite generally:

\proclaim{0.8.2 Question} Over finite fields, if $L$-equivalence
is bounded does it follow that $L$ is semi-ample?
\endproclaim

A special case of (0.8.2):

\proclaim{0.9 Question} Let $S$ be smooth surface over a finite
field. Let $D$ be a divisor on $S$. If $D \cdot C > 0$ 
for 
every irreducible curve $C \subset S$, does it follow that
$D$ is ample (or equivalently, that $D^2 > 0$)?
\endproclaim

(0.1) has a (to me) surprising consequence. 
By a {\bf map related to} $L$ I mean a proper
map $g: X \rightarrow Z$ between
algebraic spaces whose exceptional subvarieties
(i.e. the subvarieties
$W$ such that $\dim(W) > \dim g(W)$) are precisely
the $L$-exceptional subvarieties. The Stein
factorisations of any two related maps are the
same, and (if it exists) is called the {\bf map
associated to} $L$. For elementary functorial
properties of the associated map see \cite{Keel99,1.0}.
If $L$ has an associated
map we say $L$ is {\bf endowed with a map (EWM)}.
A semi-ample line bundle is EWM -- the associated
map is given by sections of a sufficiently high tensor
power, but the converse fails. However: 

\proclaim{0.10 Theorem} A nef line bundle on
a scheme projective over the algebraic closure of
a finite field is semi-ample iff it has a related
map. \endproclaim

(0.1) is a formal consequence of the existence of
certain pushout diagrams, (0.17) which though technical
is I think of independent interest. I will turn
to this next.

\subhead Content overview \endsubhead
The remainder of \S 0 is devoted to pushouts.
(0.17) is proven in \S 1. (0.1) is
proven in \S 2. The results for $\mgn$
are proven in \S 3. 
(0.10) is proven in \S 4.
An interesting example
of Koll\'ar is given in \S 5.

{\bf Thanks:} I thank James McKernan with
whom I discussed many aspects of this paper, and
who stimulated my interest by raising the
question of whether or not nef and semi-ample are
the same over a finite fields; J\'anos
Koll\'ar, who gave me a counter-example, (5.2), to
an overly
optimistic version of (0.1); Joe Harris, for
inviting me (or more precisely, accepting my
self invitation) to Harvard for a year, and
Tom Graber, with whom I jointly observed (1.9). 
I especially thank Dan Abramovich for
his very careful reading of the paper.
He found lots of mistakes, and many places
where the exposition could be simplified and
or otherwise improved.

Notations and conventions:
Throughout the note all spaces considered are algebraic
spaces of finite type over a field $k$, and all
maps are $k$-linear. I will often replace line bundles
by positive tensor powers without remark.

By a {\bf contraction} we mean a proper map
$f: X \rightarrow Y$, with $f_*(\ring X.) = \ring Y.$.

I use the following definitions and conventions of \cite{Keel99}. 
For a
line bundle $L$ on $X$ and a map $h: Y \rightarrow X$ we
write $L|_Y$ (or sometimes $L_Y$) 
for $h^*(L)$. If $L$ on $X$ is EWM then
the associated map is denoted $g_X: X \rightarrow Z_X$.

\subhead  Pushouts \endsubhead

\remark{0.12  Pushout Notation}\endremark 
We assume we are given $(X,L_X)$, with $L_X$ nef,
all spaces proper, maps $f_i: X \rightarrow X_i$,
$i = 1,2$, and line bundles $L_{X_i}$ on
$X_i$ with $L_{X_i}|_X = L_X$.

Our goal will
be to find a commutative diagram
$$
\CD
X @>{f_1}>> X_1 \\
@V{f_2}VV @V{g_1}VV \\
X_2 @>{g_2}>> P 
\endCD
$$
and a line bundle
$L_P$ on $P$ such that 
$L_P|_{X_i} = L_{X_i}^{\otimes m}$ for 
$i=1,2$ and some
$m > 0$. Furthermore we require that $L_P$ is ample,
or equivalently, by \cite{Keel99,1.0}, that $g_i$
is a related map for $L_{X_i}$, $i = 1,2$.
As a shorthand for the existence of this diagram
we will say that a {\it polarized pushout} for 
$L_{X_1},L_{X_2},f_1,f_2$ (or just $f_1,f_2$) exists.

\remark{Remarks} If a polarized pushout exists, 
then $L_X,L_{X_i}$
are semi-ample.
If $L_X$ is ample, and $f_i$ are
surjective then $f_i$ is finite, and $L_{X_i}$ is 
ample.
\endremark

By an {\bf equivalence relation} on a space $X$ 
we mean an equivalence relation on the 
$\bk$-points of $X$. We call the relation
{\bf algebraic} if it is the set of $\bk$-points
of a closed subset $R \subset X \times X$. 

We say an equivalence relation on $X$ {\bf dominates}
a map  $f: X \rightarrow Y$
if ($\bk$-points of) 
fibres of $f$ are contained in $R$-equivalence classes.
Of course if $R$ is algebraic this holds iff $R$ contains 
the reduction of $X \times_Y X$.

For a map
$h:W \rightarrow X$, and an equivalence 
relation $R$ on $X$ by the {\bf restriction $R|_W$} we mean the 
relation where $w_1$ is equivalent to $w_2$ iff
their images are $R$-equivalent. If $R$ is algebraic
this is given by the reduction of $(h \times h)^{-1}(R)$. 
For a general discussion of equivalence
relations see \cite{KeelMori97}. 

The following is immediate:
\definition{0.13 Definition-Lemma} If $f: X \rightarrow Y$
is a map and $R$ an equivalence relation on $X$, then
$R$ dominates $f$ iff $R$ is the restriction of 
an equivalence relation on $Y$. If $f$ is surjective
this equivalence relation is unique, we denote it
by $R_Y$. In this case $R_Y = (f \times f)(R)$, and
$R = R_Y|_{X}$\enddefinition

Note that for a map $h: W \rightarrow X$, the 
restriction of $L_X$-equivalence to $W$ is refined by,
but in general coarser than, $L_W$-equivalence. It's 
convenient to have a variant that is preserved by
restriction:

\definition{0.14 Definitions} We say that a closed
subset $R \subset X \times X$
is {\bf related to} a nef line bundle $L$ iff
$R$ contains any pair of $L$-equivalent $\bk$
points, and $p_2^*(L)$ is numerically
trivial on the fibres of the first projection
$p_1: R \rightarrow X$.
\enddefinition
(Note in (0.14) we do not require that $R$ is
an equivalence relation. This (slight) additional
flexibility will be useful in the proof of (3.6).)

\proclaim{0.15 Lemma} The following hold for a closed
algebraic equivalence relation $R \subset X \times X$,
a nef line bundle $L$ on $X$ and a map 
$h:W \rightarrow X$,
with $X$ and $W$ proper (over $k$). 
\roster
\item If $L$ equivalence is algebraic, then $L$
is numerically trivial on $L$-equivalence classes. 
\item If $L_X$ is semi-ample and $R$ is refined by
$L$-equivalence, 
then $R$ dominates $g_X$, and is the restriction along
the associated map
of an equivalence relation $R_{Z_X}$. If furthermore 
$R$ is related
to $L$ then either projection $R_{Z_X} \rightarrow Z_X$
is finite.
\item If $R$ is related to $L$ then $R|_W$ is related
to $L_W$.
\endroster
\endproclaim
\demo{Proof} (1) follows from \cite{BCEKPRSW,2.4}. (3)
is obvious. For (2) note that $L$-equivalence is defined by
(the reduction of) $X \times_{Z_X} X$, so obviously 
$R$ dominates $Z_X$ and so is a restriction by (0.13). 
Now suppose $R$ is related to $L$, and 
let $G$ be a connected component of an $R$-equivalence
class. Then $L|_G$ is numerically trivial, and so 
$g_X(G)$ is a point. Thus the images of $R$-equivalence
classes are finite, and the result follows, since these
are the $R_{Z_X}$ equivalence classes (and so the
fibres of $R_{Z_X}$ under either projection).
\qed \enddemo

\proclaim{0.16 Lemma} Let $R \rightarrow X$
be an algebraic equivalence relation, with
both projections $R \rightarrow X$ finite,
and $X$ reduced. Let $Z \subset R$ be a closed subset
such that both projections $Z \rightarrow X$
are generically \'etale. After replacing
$X$ by a dense $R$-invariant open set
(and $R$ and $Z$ by their restrictions), 
the equivalence relation
generated by closed points of
$Z$ is algebraic, given by a closed subset
of $R$, \'etale over $X$. 
\endproclaim
\demo{Proof}
Note, since the projections are finite, that
$R$-invariant open sets form a base for the
Zariski topology on $X$, and furthermore
finiteness of $R$ is preserved by restriction
to such sets. Replacing $Z$ 
by its union with the diagonal, and its image under the
involution $i: R \rightarrow R$ (that
switches the factors), we may assume $Z$ contains
the diagonal, and is symmetric. We inductively
define an increasing sequence of closed subsets
$Z_i \subset R$,
that are symmetric, contain the diagonal, and
are generically \'etale over $X$. 
Let $Z_1 = Z$ and let $Z'_{i+1}$ be the
reduction of 
$$
Z_i \times_{(p_2,p_1)} Z_i \subset Z_i \times Z_i.
$$
Note $Z'_{i+1}$ consists of $4$-tuples
$(a,b,b,c) \in X^4$ with $(a,b),(b,c) \in Z_i$. 
Let $Z_{i+1} \subset X \times X$ be
the reduction of the image under the (finite) 
map sending $(a,b,b,c)$ to $(a,c)$. 
It's clear that $Z_i \subset Z_{i+1}$,
and $Z_{i+1}$ is symmetric, contains the
diagonal, and is generically \'etale over $X$. 
Over generic
points of $X$ the degree of $Z_i \rightarrow X$
is bounded, thus, after shrinking $X$, the
process terminates, in an \'etale 
equivalence relation, which is by construction
the relation generated by $Z$. \qed \enddemo

The main technical observation of this 
note is the following:

\proclaim{0.17 Theorem: Polarized pushouts over finite fields} 
If the base field $k$ is the algebraic closure of a finite
field, then given $L_{X_i},f_i$ as in (0.12), with
$L_{X_i}$ semi-ample, a polarized pushout exists iff
there is an algebraic equivalence relation $R$ 
related to $L$ and dominating $f_1,f_2$. Furthermore,
in this case there exists a pushout with 
$X \rightarrow P$  dominated by $R$.

\endproclaim

To understand the proof of (0.17) the reader may
wish to look first at the proof of the
special case (1.9) where
the role of the assumptions, especially in light
of example (5.2), becomes clear. (0.17) is reduced
to (1.9) by a series of formal manipulations of the 
sort used in \cite{Keel99,\S 2}, and
\cite{Koll\'ar97}.

\subhead \S 1 Proof of (0.17) \endsubhead

\proclaim{1.1 Lemma} To prove (0.17) in general, it is
enough to prove it in the case when $L_X$
and $L_{X_i}$ are all ample.
\endproclaim
\demo{Proof} Suppose (0.17) holds in the ample case.
Consider the general case. Note $g_{X_i} \circ f_i$ 
is related to $L_X$ and so factors through $g_X$. Thus
we have a pair of finite maps $f_i': Z_X \rightarrow Z_{X_i}$ 
with $f_i' \circ g_X = g_{X_i} \circ f_i$. By (0.15)
$R$ is the restriction of some $R_{Z_X}$ on $Z_X$,
easily seen to dominate $f_i'$.
$L_X$, $L_{X_i}$ are pullbacks of ample line
bundles $L_{Z_X}$, $L_{Z_{X_i}}$. The assumptions
of (0.17) are thus satisfied 
for the induced diagram
$$
\CD
Z_{X} @>{f_1'}>> Z_{X_1} \\
@V{f_2'}VV    @. \\
Z_{X_2}. @.  @. 
\endCD
$$
So by the ample case we 
have a polarized pushout 
$$
\CD
Z_X @>{f_1'}>> Z_{X_1} \\
@V{f_2'}VV @V{g_1'}VV \\
Z_{X_2} @>{g_2'}>> P 
\endCD
$$
with the line bundles pulled back from ample $L_P$ on $P$,
and $R_{Z_X}$ dominating $P$. As $R$ is restricted from
$Z_X$ it dominates $P$ as well, so 
$g_i' \circ g_{X_i}$ give the required polarized pushout.
\qed \enddemo

\proclaim{1.3 Lemma} Let $f: Y \rightarrow X$ be 
a finite map, with $Y$ quasi-projective (and
$X$ an algebraic space). Then the
natural map 
$H^1(f_*(\uring Y.)) \rightarrow \pic Y.$ is an
isomorphism (where $H^1$ means \'etale
cohomology).
\endproclaim
\demo{Proof} By the Leray spectral sequence
it's enough
to show for each $x \in X$ the
restriction to the stalk
$$
\pic Y. \rightarrow R^1f_*(\uring Y.)_x
$$
is trivial. For this we may replace $X$ by an
\'etale neighborhood of $x$ and so may assume
$X$, and hence $Y$, are affine. But now given
any line bundle $A$ we can find a section 
$\sigma$ non-vanishing at any point of (the finite
set)
$f^{-1}(x)$. Let $Z$ be the zero locus of $\sigma$.
$$
x \in U := f(Z)^c \subset X
$$
and $A|_{f^{-1}(U)}$ is trivial. \qed \enddemo

\proclaim{1.4 Lemma} If $f_1$
is a closed embedding, $f_2$ is 
finite and the $L_{X_i}$ 
are ample
then a polarized pushout $P$ exists,
with $g_2$ a closed embedding and $g_1$ finite.
The pushout diagram is a pullback (i.e. 
$X = g_1^{-1}(X_2)$), and outside of $X$,
$g_1$ is an isomorphism onto its image.
Any equivalence relation dominating $f_2$ 
will dominate $P$. \endproclaim
\demo{Proof}  The final claim follows from the claim 
that precedes it.
Replacing $X_2$ by the scheme-theoretic
image of $f_2$, and running the argument twice
(the second time with each $f_i$ a closed 
embedding) we can assume $f_2$ is surjective. 
By \cite{Artin70,6.1} a universal pushout
exists with the required properties 
in the category of algebraic
spaces. 
Thus it's enough to show that some positive tensor 
power of
$L_{X_i}$
descends to $P$.
Consider the complex of \'etale sheaves
of Abelian groups
$$
1 \rightarrow \uring {P}. \rightarrow 
{g_2}_*(\uring X_2.)\times 
{g_1}_*(\uring X_1.) \rightarrow
({g_i} \circ f_i)_* (\uring X.) \rightarrow 1.
$$
By (1.3) and the long exact cohomology sequence
it is enough to show this is exact.
At the left and middle this
follows from the universal property of
$P$ (and its restriction to \'etale open
sets) applied to maps to ${\Bbb A}^1$,
see e.g. \cite{Koll\'ar97,8.1.3}. For exactness
on the right, it's enough to show 
$$
{g_1}_*(\uring X_1.) \rightarrow (g_1 \circ f_1)_*
(\ring X.)
$$
is surjective. For this we can assume $P$, and
hence all the spaces, are affine. A unit on 
$X$ lifts to a function on $X_1$, invertible
in a neighborhood of $X$. Since 
$X = g_1^{-1}(X_2)$, this neighborhood is
the inverse image of a neighborhood of
$X_2$ in $P$. \qed \enddemo

\proclaim{1.5 Lemma} Assume positive
characteristic. If 
$f_1$ is a finite universal homeomorphism, and
the line bundles  $L_X$, $L_{X_i}$
are all ample, then 
a polarized 
pushout exists, with $g_2$ a finite universal 
homeomorphism.
Any equivalence relation dominating $f_2$ will
dominate $P$.
\endproclaim
\demo{Proof} By \cite{Koll\'ar97,8.4} the pushout exists,
and line bundles descend by \cite{Keel99,1.4}. The
condition that an equivalence relation dominate
$P$ is set-theoretical, so since $P$ is homeomorphic
to $x_2$, the final claim is clear. \qed \enddemo

\proclaim{1.6 Lemma} Assume positive
characteristic. Suppose there is a finite universal
homeomorphism $h: X' \rightarrow X$ such that a polarized
pushout for $h \circ f_i,L_X|_{X'}$ exists. Then a polarized
pushout for $f_i,L_{X_i}$ exists. \endproclaim
\demo{Proof} Immediate from \cite{Keel99,2.1} and
\cite{Keel99,1.4}. \qed \enddemo

\proclaim{1.7 Lemma} Assume $k$ is the algebraic
closure of a finite field. Let $f:Y \rightarrow X$ be 
a proper surjection between projective schemes. The
kernel of $f^*:\pic X. \rightarrow \pic Y.$ is torsion.
\endproclaim
\demo{Proof} Any line bundle in the kernel will be numerically
trivial and thus torsion, see e.g. \cite{Keel99,2.16}. 
\qed \enddemo

\proclaim{1.8 Lemma} Assume $k$ is the algebraic
closure of a finite field. 
Suppose there is a proper surjection 
$h: X' \rightarrow X$ such that a polarized
pushout for $h \circ f_i,L|_{X'}$ exists. Then a polarized
pushout for $f_i,L_{X_i}$ exists. Furthermore if there is
an equivalence relation
$R$ on $X$ dominating $f_i$, and there is a pushout
for $h \circ f_i$ dominated by $R|_{X'}$ then there
is a pushout for $f_i,L_{X_i}$ dominated by $R$.
\endproclaim
\demo{Proof} Assume we have a polarized pushout
$g_1,g_2$ for $h \circ f_i$. To obtain the desired polarized
pushout we need only show
that $g_1 \circ f_1 = g_2 \circ f_2$, and $L_P|_X = L_X$.
By (1.7) and \cite{Keel99,2.1} these are achieved
after replacing line bundles by powers, and $P$ by
a finite universal homeomorphism $P \rightarrow P'$. Now
the final remarks with respect to $R$ follow (as we have
only modified the spaces by universal homeomorphism and
the conditions on the equivalence relations are set-theoretical).
\qed \enddemo

\proclaim{1.9 Lemma} (0.17) holds over any field 
assuming $X,X_i$
are normal, $L_{X_i}$ are ample, and $f_i$ are
finite surjections.
\endproclaim
\demo{Proof} We are free to replace $X$ by
a finite cover $h:X' \rightarrow X$: 
Since $X$ is normal
the kernel of the pullback on Picard groups is torsion,
and since $X$ is reduced 
$$
\operatorname{Hom}(X,T) @>{h \circ}>> 
\operatorname{Hom}(X',T)
$$
is injective for any $T$, so we can argue as in the proof
of (1.8). 

Note since $R$ is related to an ample line bundle, 
either projection $R \rightarrow X$ is finite.

By (1.5), factoring the field extensions
into separable and purely inseparable parts, we
may assume the field extensions $K(X_i) \subset K(X)$
are separable. 

Let $K$ be the intersection of function fields
$K(X_1) \cap K(X_2) \subset K(X)$. Suppose
for the moment that $K \subset K(X)$
is finite. 
Let $X'$
the integral closure of $X$ in the Galois closure
of $K(X)$ over $K$. We may replace $X$ by
$X'$, thus we may assume $K(X_i) \subset K(X)$ are Galois
extensions, and so $f_i$ is the quotient by
the Galois group, $G_i$. 
Let $H$ be the subgroup of $\Aut(X)$ generated
by $G_1$ and $G_2$. 
$H$ is a subgroup of $\Aut(K(X),K)$ and so finite.
We can take the geometric quotient $P = X/H$. 
Some power of $L_X$ will descend to $P$. $R$ will
contain the equivalence relation generated by
$G_1,G_2$ which is exactly (the reduction of)
$X \times_P X$.

Thus it is enough to show that $K \subset K(X)$
is a finite extension.
Since $R$ dominates $f_i$ any
$R$-invariant open
set is the inverse image of its image under $f_i$, so
to check the finiteness of $K \subset K(X)$ we can
replace $X$ by any non-empty $R$-invariant set, and
such sets form a base for the topology, e.g. as 
either projection $R \rightarrow X$ is finite. Note
also that the finiteness of $R \rightarrow X$ is 
preserved by
restriction to $R$-invariant open sets. 
Let $Z \subset R$ be the reduction of 
the union of 
$X\times_{X_1} X$ and $X \times_{X_2} X$.
Either projection $Z \rightarrow X$ is generically
\'etale, since the extensions 
$K(X_i) \subset K(X)$ are separable, so
by (0.16) we may assume 
there is an algebraic 
equivalence relation $R' \subset R$ 
dominating $f_1,f_2$, with $R' \rightarrow X$
finite and \'etale. 
A geometric quotient by $R'$ exists 
(in fact it's just the algebraic
space defined by $R'$), 
see e.g. \cite{KeelMori97,4.8,5.1}, a 
finite map $q:X \rightarrow X/R'$.
Since $R'$ dominates $f_i$, $q$ factors through $f_i$ and
so $K(X/R') \subset K$. Since $K(X/R') \subset K(X)$
is a finite extension, so is $K \subset K(X)$.
\qed \enddemo

\demo{Proof of (0.17)} We prove the result by induction 
on the dimension of $X$. By (1.1) we may assume all
of $L_X$, $L_{X_i}$ are ample.
By (1.8) we can assume $X$ is normal.
By (1.4) applied to scheme-theoretic images, we may
assume the $f_i$ are surjective and $X_i$ reduced.
By (1.9) we 
reduce to the case when $f_1$ is the normalization:
Indeed by (1.9) we have a diagram
$$
\CD
X @>{{\til f.}_1}>>  \til X_1. @>{p_1}>> X_1 \\
@V{{\til f.}_2}VV   @V{{\til g.}_1}VV  @. \\
\til X_2.  @>{{\til g.}_2}>> P @. \\
@V{p_2}VV @. @. \\
X_2 @. @. @. @.
\endCD
$$
with $p_i$ the normalization, where the square is
a polarized pushout relative to $R$.  
By (0.15), $R = R_{\tX_i}|_X$
for a finite equivalence relation on ${\til X_i.}$ dominating
$p_i,{\til g.}_i$. So by (0.15), it's enough
to construct a polarized pushout $h: P \rightarrow P'$
of $p_1,{\til g.}_1$, dominated by $R_{\tX_1}$,
and then a polarized pushout of 
$p_2,h \circ {\til g.}_2$, dominated by $R_{\tX_2}$.

We write $p = f_1$. 
Recall that the conductor is the sheaf
of functions on $X$ which multiply
$\ring X.$ into $\ring X_1.$. This is at once an
ideal of $\ring X.$ and of $\ring X_1.$. Let
$C \subset X$, $D \subset X_1$ be the associated
subspaces. The definitions imply that the
diagram
$$
\CD
C @>i>>  X  \\
@VpVV   @V{p}VV \\
D @>i>> X_1
\endCD
\tag{1.10}
$$
(where the maps $i$ are the closed embeddings
defined by the conductor ideals)
is a fibre diagram, and a universal pushout, see e.g.
\cite{Reid94,2.1.}.
By induction on dimension we have a polarized pushout 
diagram 
$$
\CD
C @>{f_2 \circ i}>> f_2(C)\\
@V{p}VV   @VVV \\
D @>>> P'
\endCD
$$
for $L_{X_2},L_{X_1}|_D$, dominated by $R|_C$. Further
by (1.4) we have a polarized pushout 
$$
\CD
f_2(C) @>{i}>> X_2 \\
@V{p}VV  @V{p}VV \\
P' @>{i}>> P 
\endCD \tag{*}
$$
with $i$ closed embeddings, $p$ finite surjections,
which is also a pullback diagram, with $X_2 \rightarrow P$
an isomorphism outside of $f_2(C)$.
In particular we have a commutative diagram of finite maps
$$
\CD
C @>{f_2 \circ i}>> X_2 \\
@VpVV   @VVV \\
D @>>> P
\endCD
$$
and compatible ample $L_P$. 
Since the conductor diagram is a
universal pushout, the composition 
$$
X \rightarrow X_2 \rightarrow P
$$
factors through $p:X \rightarrow X_1$. $L_{X_1}$
descends to $L_P$ (after taking powers) by (1.7). Finally
we check this pushout is dominated by $R$. Suppose
$x_1,x_2 \in X$ are distinct points that 
have the same image in $P$. We show they are $R$-equivalent.
We may assume $f_2(x_1),f_2(x_2)$ are distinct (or the
$x_i$ are $R$-equivalent by assumption). 
Since $p: X_2 \rightarrow P$ is an isomorphism
outside of $f_2(C)$, and (*) is a pullback, we have
$f_2(x_i) \in f_2(C)$. So replacing $x_i$ by 
$X \times_{X_2} X$, and hence $R$, equivalent points,
we can assume $x_i \in C$ and so $R|_C$ equivalent.
\qed \enddemo

\subhead \S 2 Proof of (0.1) \endsubhead
The following is clear from the definitions:
\proclaim{2.1 Lemma} Let $W \rightarrow X$ 
be a map between proper algebraic spaces.
Let $L$ be a nef line bundle on $X$.
$L_W$ equivalence is a refinement of
the restriction of $L_X$-equivalence. 
\endproclaim

\proclaim{2.2 Theorem} 
Let $p:\tX \rightarrow X$ be
a surjection between algebraic spaces
proper over 
the algebraic closure of a finite field. Let $L$ be a 
line bundle on $X$. Then $L$ is semi-ample iff 
the following two conditions hold:
\roster
\item $p^*(L)$ is semi-ample, and 
\item there is a closed
algebraic equivalence relation $R \subset X \times X$
related to $L$ (i.e. $R$ is refined 
by $L$-equivalence, and $L$ is 
is numerically trivial on $R$-equivalence classes).
\endroster
\endproclaim
\demo{Proof} If $L$ is semi-ample then clearly
the conditions must hold -- for $R$ one takes
$X \times_{Z_X} X$. We prove the other implication by
induction on the dimension of $X$. 
Since we may replace $\tX$ by anything that surjects
properly onto it, we can assume $\tX$ is normal. 
By \cite{Keel99,1.4}
(passing to the reduction) we may assume $X$ is 
reduced. (0.1) holds for $X$ normal, over any
field, and without the boundedness assumption,
see e.g. the proof of \cite{Keel99.2.10}.
(Alternatively, in positive characteristic, one can apply
\cite{Keel99,2.10} to the case where the 
base space is normal and $D$ and $C$ are empty.)
So the restriction of $L$ to the normalisation
of $X$ is semi-ample, and thus we may assume 
from the start that $X$ is reduced and 
$\tX$ is the normalisation. 

We let $C,D$ be
the conductors as in (1.10) (but here 
$p:\tX \rightarrow X$ replaces $p:X \rightarrow X_1$).
By (2.1) $R|_D$ is refined by $L_D$-equivalence, and
$L_D$ is clearly numerically  trivial on $R|_D$ equivalence
classes. Thus
by assumption and induction, 
$L_{\tX},L_{C},L_{D}$ are
all semi-ample.  Any two points in the same fibre of
the composition
$$
g_{\tX} \circ i: C \rightarrow \tX \rightarrow Z_{\tX}
$$
map to $L_X$-equivalent points of $X$, and thus are
$R|_C$ equivalent, i.e. $R|_C$ dominates $Z_{\tX}$.
$R|_C$ obviously dominates $C \rightarrow D$
(since $R|_C$ is restricted from $D$).
Thus by (0.17) there is a polarized pushout 
$$
\CD
C @>>> Z_{\tX} \\
@VVV   @VqVV \\
D @>>> P
\endCD
$$
with $L_{Z_{\tX}}$ (and thus $L_{\tX}$) pulled back
from an ample line bundle $L_P$. Since the conductor
diagram  is a universal
pushout, $q \circ g_{\tX}$ factors through $X$.
By (1.7), as a $\bbbq$-line bundle $L_X$ is the
pullback of $L_P$, thus semi-ample. \qed \enddemo

\proclaim{2.3 Lemma} Let $p: W \rightarrow X$
be a surjection between proper algebraic spaces.
If $p^*(L)$-equivalence
is algebraic and $L$-equivalence is
bounded, then $L$-equivalence is algebraic. \endproclaim
\demo{Proof} Let $m$ be the bound for $X$. Let
$S' \subset W \times W$ define $L_W$ equivalence,
and let $S \subset X \times X$ be its image. 
Clearly the $k$-points of $S$ are contained in
the union of $L$-equivalent pairs. 
Note that if $x_1,x_2 \in C$ and $C \cdot L =0$
for irreducible $C$, then $(x_1,x_2) \in S$,
since we can lift $C$ to an $L_W$-trivial 
irreducible curve.
So by assumption, $x,y$ are 
$L$-equivalent iff we can find
$x = x_1,x_2,\dots,x_m = y$ with $(x_i,x_{i+1}) \in S$. 
It's clear this defines a closed subset of $X \times X$. \qed
\enddemo
\demo{Proof of (0.1)} Immediate from (2.2 - 2.3). \qed \enddemo

\proclaim{2.4 Lemma} Let $L$ be a nef line bundle
on a complete algebraic space $X$.
If there is a closed subset $W \subset X \times X$
related to $L$, then $L$-equivalence is bounded.
\endproclaim
\demo{Proof} We can take for the bound $m$ the maximum number
of irreducible components in any fibre of the
first projection $W \rightarrow X$. \qed \enddemo

\subhead \S 3 Semi-ampleness on $\mgn$ \endsubhead

{\bf Assumptions Global to \S 3} 
Throughout this section $L$ is a nef line bundle on
$\mgn$. $k$ is a fixed finite field. All spaces
are assumed to be of finite type over $k$ and
all morphisms are $k$-linear.

{\bf 3.0 Topological Stratification and Product Decompositions:}
$\mgn$ has a tautological stratification by
topological type. Here we recall the main
points, see \cite{Keel99,pg. 274} for
more precise details. 
A codim $i$-stratum is a connected
component of the locus corresponding to stable pointed
curves with precisely $i$-singular points. I will
abuse language and refer to the closure of a stratum
as a {\bf closed stratum}. 
There is one stratum for each topological
type of stable pointed curve, or equivalently, for
each isomorphism class of decorated dual curve --
where by decorated I mean the vertices are labeled by
genera, and labeled points are indicated by labeled
edges, incident to a unique vertex (of course corresponding
to the component of the curve containing the labeled point).
For
a stable pointed curve $[E] \in \mgn$, we
write $X_E$ for the unique closed stratum containing
$[E]$ in its interior. Of course this depends
only on the topological type of $E$.
We write $X_E^{\circ}$ for the (open) stratum.

The normalisation of each closed stratum is a natural
projective moduli space. More convenient for our purposes is
a certain branched cover which
we call the {\bf product decomposition} of the stratum:
Let
$E$ be a stable
pointed curve, and $X = X_E$. 
Let 
$\Gamma$ be the dual graph of $E$. Let
$P_X$ be the product of 
$\vmgn g_i. N_i.$, one for each vertex, with
$g_i$ the genus of the vertex, and $N_i$ the
(ordered) set of incident edges. In the notation of
\cite{Keel99,pg. 274}, 
$$
\align
P_X &= \overline{\Cal M}^{X \cup Y}_{E} \\
N_i &= X_i \cup Y_i.
\endalign
$$
Note formation of $P_X$ requires ordering the vertices of
$\Gamma$, and then for each vertex the collection
of unlabeled incident edges. In this sense it is not
unique, though any two choices are isomorphic.
$H=\Aut(\Gamma)$ (meaning automorphisms of the
decorated graph) 
acts naturally on $P_X$ and the quotient
is $\tX$, the normalisation of $X$. See 
e.g. 
\cite{HainLooijenga97,\S 4}, or 
\cite{GraberPandharipande, \S A}.
The quotient
of the interior of $P_X$ is $X^o$ (where
the interior of $P_X$ means the product
of the interiors of its factors).

We note that
$\mgn$, all of its $\bQ$-line bundles,
the stratification by topological type, the
closed strata, their normalisations, 
and product decompositions of all strata, are all
defined over $k$ (indeed over $\bZ$), see 
\cite{DeligneMumford69}, \cite{Moriwaki01}.

In various notations we replace
$X_E$ by $E$, e.g. we will sometimes write
$P_E$ for the product decomposition $P_{X_E}$.

{\bf 3.0.1:} 
Suppose that $F$ is a degeneration of the stable pointed
curve $E$, so $X_F \subset X_E$. Observe that
$P_F \rightarrow X_F \subset X_E$ factors
(non canonically) through 
$P_E \rightarrow X_E$: Let $E_i$ be the irreducible
components of $E$. 
$F$ can be written as a union of connected
(not necessarily irreducible)  curves $F_i$, with
$F_i$ a degeneration of $E_i$. Let $\tE_i$
be the normalisation of $E_i$. There is 
a corresponding connected partial normalisation $\tF_i$ 
of $F_i$, $\tF_i$ a degeneration of $\tE_i$.
Each $E_i$ corresponds
to a factor, $\vmgn g_i. N_i.$, of $P_E$. 
$\tF_i$ gives a closed stratum $X_{\tF_i}$ of $\vmgn g_i. N_i.$.
The composition
$$
\underset {i} \to \times P_{\tF_i} \rightarrow \underset {i} \to \times \vmgn 
g_i. N_i. =P_E
\rightarrow X_E
$$ 
has image $X_F$. $\times_i P_{\tF_i} \rightarrow X_F$
is the product decomposition $P_F$. 

The goal here is to give a necessary and sufficient
condition for $L$ to be semi-ample. We will use (0.1),
and (2.4),
and so look for an algebraic equivalence relation $R$ 
related to $L$  (see (0.14) for a reminder of
what it means for a closed subset to be related
to a nef line bundle). 
As the argument is a bit technical,
let me begin by sketching the general philosophy: 

If one can guess the set theoretic 
fibres of the associated map, then one can write
down (a candidate) $R$ directly (note we only
need $R$ related to $L$-equivalence, so its
equivalence classes need not be precisely the
fibres of the associated map). 
For example, for $\lambda$, the
associated map is to the Satake compactification
of the moduli of Abelian varieties,
so we say $[C],[D]$ are $R$-equivalent iff
the union of the non-rational
components 
of their respective normalisations are isomorphic. 
For $\psi_i$: For a pointed curve $C$ let $C_i$ be the
irreducible component containing the $i^{th}$ point, and
let $\tC_i,S$ be its normalisation, together with the
(unordered) set $S$ of all
its distinguished points, i.e. all the labeled
points lying on $C_i$, together with the inverse images
of all the singular points of $C$.  Then we say
$[C],[D]$ are $R$-equivalent iff 
$(\tC_i,S)$ and $(\tD_i,S)$ are isomorphic.
(In fact (a slight variant) of this 
equivalence relation was introduced purely
topologically (in characteristic $0$) by Kontsevich.
Looijenga observed that it was algebraic, and in
fact defined $\psi_i$ equivalence, see \cite{Looijenga95}.) 
In these examples one can
argue $R$ is algebraic by using the associated
moduli spaces. E.g. for $L = \psi_i$, let
$\cQ$ be the disjoint union of $\vmgn g. n./S_n$,
one copy for each $g,n$,
and $\cX^{o}$ the disjoint
union of the (open) strata. We have a natural
map $\cX^o \rightarrow \cQ$, sending $[C]$
to $[\tC_i,S]$ (notation as above). Since
$\cQ$ is the moduli space of stable unordered
pointed curves, two points of $\mgn$
are (by definition) $R$-equivalent
iff they have the same images in $\cQ$, thus 
$R$ is the set of closed points of 
$\cX^o \times_{\cQ} \cX^o \subset \mgn \times \mgn$.
(Of course it remains to show that $R$
is closed and related to $\psi_i$, this requires
some work, see the proof of (3.6).)

For general $L$ we do something similar. 
We will know by induction on dimension that
$L$ is semi-ample on the normalisation of any
proper closed stratum. Then we can form a space
analogous to $\cQ$ above out of the images
of the associated maps. To do so we
make use of a simple construction:

By a polarised variety we mean a pair $(Z,L_Z)$ of
a projective $Z$ and an ample $\bbbq$-line bundle
$L_Z$. Note since $k$ is finite, a $\bbbq$-line
bundle is the same as its Chern class in the
Neron-Severi group with rational coefficients.

\proclaim{3.1 Lemma} Let $(Z,L_Z)$ be 
a polarized variety (projective over $k$). 
The discrete group $G =  \Aut(Z,L_Z)(k)$ is
finite. Let $Q_Z = Z/G$. $L_Z$
descends (as a $\bbbq$-line bundle) to an ample
$L_{Q_Z}$.  
Let $f: (Z,L_Z) \rightarrow (W,L_W)$ be an
isomorphism of polarised varieties. 
Then there is an
induced isomorphism of polarised
varieties $Q_Z \rightarrow Q_W$, which
is independent of the choice of $f$. \endproclaim
\demo{Proof} $\Aut(Z,L_Z)$ is of finite type
(for any polarized variety), so its set of 
rational points over the finite field $k$ is finite.
The rest is easy. \qed \enddemo

\definition{3.2 Definition} Let $X$ be
a closed stratum such that 
$L_{P_X}$ is semi-ample. 
We define
$$
(Q_X,L_{Q_X}) := (Q_{Z_{P_X}},L_{Q_{Z_{P_X}}})
$$
where the right hand side uses the notation of (3.1). 
\enddefinition

\proclaim{3.3 Lemma} Notation as in (3.2). The
composition
$$
q_{P_X}: P_X \rightarrow Z_{P_X} \rightarrow Q_X
$$
factors through $P_X \rightarrow \tX$. In
particular there 
is an induced map $q_X: X^o\rightarrow Q_X$.
$q$ and $(Q_X,L_{Q_X})$ depend only on $X$ and $L$.
\endproclaim
\demo{Proof} Let $X = X_E$, for a topological
type of stable pointed curve $E$. 
Let $H$ be, as in (3.0), the automorphisms
of the dual graph of $E$. Since $L_{P_X}$
is pulled back from $X$ (and so from $\tX$)
it has a natural $H$ linearization. Thus
$H$ acts on 
$(Z_{P_X},L_{Z_{P_X}})$. The quotient
is $Z_{\tX}$. Indeed:
$$
\align
Z_{P_X} / H & = \proj(R(P_X,L_{P_X})^H) \\
            &= \proj(R(\tX,L_{\tX})) \\
            &= Z_{\tX} \endalign
$$
(where $R$ indicates the graded ring of sections). 
Thus $Z_{P_X} \rightarrow Q_{X}$ factors
through $Z_{\tX}$, and thus 
$P_X \rightarrow Q_X$ factors through $\tX$, inducing
$q: \tX \rightarrow Q_X$. The 
choices here were the orderings (of factors,
and marked points) in the construction of $P_X$.
Different choices yield isomorphic polarized
varieties and so, by (3.1), the same $Q_X$. \qed \enddemo

Now in the proof of (3.6) we will define $[E]$ and $[F]$
to be $R$-equivalent (assuming $L_{P_E}$ and
$L_{P_F}$ are semi-ample) iff  $Z_{P_F}$ and $Z_{P_E}$ are
isomorphic polarized varieties, and $q_E([E])$ and $q_F([F])$,
($q_E = q_{X_E}$ of (3.3)) are identified under the
canonical identification $Q_E = Q_F$ of (3.1-2).
To show $R$ is closed and related to $L$ we 
will need
some further conditions, which we formalize as
follows:

\remark{3.4 Compatibility Condition} Let $E$ be
the topological type of a stable pointed curve.
We say (3.4) holds for $E$ if there exists a 
contraction 
$p: P_E \rightarrow W$ with the following
properties:
\roster
\item $L_{P_E}$ is pulled back from $W$.
\item $p$ contracts any irreducible 
$L_{P_E}$-exceptional curve
which meets the interior of $P_E$ 
\item Assume $L_{P_E}$ is semi-ample. 
Let $[F] \in X_E$ , 
and let $[\tF],[\tE] \in P_E$ be
points with image $[F],[E] \in X_E$. 
If $p([\tF]) = p([\tE])$ then the polarized
varieties $Z_{P_E}$ and $Z_{P_F}$ are
isomorphic and $q_F([F]) = q_E([E])$ in
$Q_F = Q_E$, under the identification
(3.1) and the maps, $q_E,q_F$, of (3.3). 
\endroster
\endremark

Remark: Note in (3.4.3) that since $[F] \in P_E$,
there is, by (3.0.1), a factorisation $P_F \rightarrow P_E$,
and thus $L_{P_F}$ is semi-ample as well, so 
$q_F: X^0_F \rightarrow Q_F$ is defined.

\proclaim{3.5 Lemma} Let $[E] \in \mgn$ be a stable pointed
curve. Suppose for each 
$h: \vmgn r. N. \rightarrow \mgn$ coming from
a factor of the product decomposition $P_E \rightarrow X_E$,
either $h^*(L)$ is numerically trivial, or 
there is 
a subset $S \subset N$ such that 
$$
h^*(L) = \pi_S^*(L_S)
$$
for a nef line bundle $L_S$ on $\vmgn r. S.$ with
$\ex L_S. \subset \partial \vmgn r. S.$
(here we require that $S$ is such that $\vmgn r. S.$
is defined, i.e. if $r =1$, $S$ is non-empty,
and if $r =0$, $S$ has at least $3$ elements).

Then (3.4) holds for $E$.
\endproclaim
\demo{Proof} 
If $L_{P_E}$ is numerically trivial we take
$W$ to be a point and there is 
nothing to check (the spaces $Q_F$ and $Q_E$
in (3.4.3) will both be points). So we
assume otherwise.

Let $W$ be 
the product of the $\vmgn r. S.$, over those
factors of the product decomposition where
$p^*(L)$ is not numerically trivial.
Let  $p: P_E \rightarrow W$ be the product of (projection
to the factor composed with) the
maps $\pi_S$. 
By \cite{GKM00,1.1},
$L_{P_E}$ is a tensor product of (the
pullbacks of) the $h^*(L)$ (as $h$ varies over
the factors of $P_E$). Thus (3.4.1)
holds. 

Let $C$ be an irreducible $L_{P_E}$ exceptional
curve that meets the interior of $P_E$.  Note 
$$
\pi_S^{-1}(\partial \vmgn r. S.) \subset
\partial \vmgn r. N. .
$$
Thus $\pi(C)$ meets the interior of $\vmgn r. S.$.
By
the projection formula $c_1(L_S) \cdot \pi_S(C) =0$
for all $S$. So by the assumption 
$\ex L_S. \subset \partial \vmgn r. S.$,
$\pi_S(C)$ is a point for all $S$, whence (3.4.2).

Now suppose $L_{P_E}$, or equivalently,
$L_W$, is semi-ample, 
and $\pi_S([\tF]) = \pi_S([\tE])$, for
$F$ a degeneration of $E$. Choose
a factorisation $P_F \rightarrow P_E$ as
in (3.0.1). Assume for the moment that 
the composition $P_F \rightarrow W$ is 
a contraction:

Then since
$L_{P_{F}}$ and $L_{P_E}$ are pulled back from $L_W$ 
(the tensor product of the $L_S$)
along contractions, the images of the associated
maps, $Z_{P_F}$, $Z_{P_E}$, and $Z_W$,  are 
isomorphic polarised varieties. Thus $Q_F$ and
$Q_E$ are canonically identified by (3.1-2).
By definition (3.3) the images $q_E([E])$ and
$q_F([F])$ in $Q_F = Q_E$ are $q_{P_F}([\tF])$
and $q_{P_E}([\tE])$. $q_{P_F}$ and $q_{P_E}$
factor through $p: P_E \rightarrow W$, and 
$p([\tE]) = p([\tF])$ by assumption, whence
(3.4.3).

So
it's enough to show $P_F \rightarrow W$ 
is a contraction. We first reduce to the
case when $E$ is irreducible. 
We follow the notation of
(3.0.1), so $P_E = \times_i \vmgn g_i. N_i.$.
For factors with $L|_{\vmgn g_i. N_i.}$ 
non-trivial, let $S_i \subset N_i$ be the 
corresponding subset (called $S$ in the statement
of (3.5)). Let $W_i = \vmgn g_i. S_i.$
and let $p_i:P_{\tF_i} \rightarrow W_i$ 
be the composition of 
$$
P_{\tF_i} \rightarrow X_{\tF_i} \subset \vmgn g_i. N_i.
$$
with 
$$
\pi_{S_i}: \vmgn g_i. N_i. \rightarrow \vmgn g_i. S_i..
$$
$P_F = \times_{i} P_{\tF_i}$ and 
$P_F \rightarrow W$ is the product of the $p_i$. So
it's enough to show each $p_i$ is a contraction. But
note if from the start we replace $\mgn$ by
$\vmgn g_i. N_i.$ and $E,F$ by $\tE_i,\tF_i$, then
our construction replaces $P_F \rightarrow W$ by $p_i$. So we
may assume from the start that $E$ is irreducible.

Now $P_E$ has a single
factor $\vmgn r. N.$, $p = \pi_S$ and
$W = \vmgn r. S.$.

$\pi_S$ sends $[C,N]$ to the stabilisation
of $[C,S]$. Note since $E$ is irreducible,
$[\tE,S]$ is stable, and by assumption,
is the stabilisation of $[\tF,S]$.
So by the definition of 
stabilisation 
there  is a unique irreducible
component $\tF_S$ of $\tF$ 
such that: each connected component of 
$\tF \setminus \tF_S$ meets $S$ in at most one point,
$\tF_S$ is smooth, $\tF \setminus \tF_S$ is a 
(possibly disconnected) tree of smooth
rational curves, and $\tF_S$ is the underlying curve of $\tE$.
It follows that the closed stratum 
$X_{\tF} \subset \vmgn r. N.$ is
normal, equal to its product decomposition, $P_F$,
and is a product of 
$\vmgn r. S \cup T.$ with a product of various
$\vmgn 0. N_i.$, where $T$ are the labels 
in $N \setminus S$ that lie 
on $\tF_S$, and under this identification
$X_{\tF} \rightarrow \vmgn r.S.$ is a 
contraction, the 
composition of projection onto the
$\vmgn r. S \cup T.$ factor with the natural 
contraction
$\vmgn r. S \cup T. \rightarrow \vmgn r. S.$. 
\qed \enddemo

The next proposition is a restatement of (0.3.1) and
(0.4): 

\proclaim{3.6 Proposition} Let $L$ be a nef line bundle
on $\mgn$. Assume the characteristic
$p > 2$. $L$ is semi-ample provided one of the
following holds:
\roster
\item $L$-equivalence
is bounded and $r^*(L)$ is semi-ample for
$$
r: \vmgn 0. 2g+n. \rightarrow \mgn
$$
the natural map given by the rational locus, or
\item For each $h: \vmgn r. N. \rightarrow \mgn$
a component of a product decomposition of a 
closed stratum,
with $r \leq 1$, either
$h^*(L)$ is trivial, or 
$h^*(L) = \pi_S^*(L_S)$ for
some subset $S \subset N$ and some 
nef line bundle $L_S$ with 
$\ex L_S. \subset \partial \vmgn r. S.$.
\endroster
\endproclaim
Remark: In (3.6.2) we require, as in
(3.5) that $S$ is such that $\vmgn r. S.$
is defined. We allow the possibility 
that $S = N$, in which case the assumption is that
$\ex {h^*(L)}. \subset \partial \vmgn r. N.$. 
\demo{Proof} 
We induct on the dimension of $\mgn$. It's clear
that the assumptions apply to $h^*(L)$ for any 
component of a product decomposition
of a proper stratum. 
$L_{P_X}$ is a tensor product
of nef line bundles pulled back from
the factors, see \cite{GKM00,1.1}. Thus $L_{P_X}$ is by
induction semi-ample for any proper closed 
stratum. In particular the restriction of
$L$ to the normalisation of $\partial \mgn$ is
semi-ample. 

We can assume that $L$ is big and
$\ex L. \subset \partial \mgn$: Suppose
not. If $g \geq 1$ then by 
\cite{GKM00,0.9}, $L$ is pulled
back from some lower dimensional 
$\vmgn g. S.$, and we can apply
induction. If $g =0$, case (1) is vacuous. 
In case (2), by the assumption 
applied to the entire space,
$L$ is again pulled back.

Thus the restriction of $L$ to the
normalisation of $\ex L.$ is semi-ample, so 
by (0.3) it's enough to show that $L$-equivalence
is bounded. This completes the proof for (1), so we 
can assume that we have (2). Let $E$ be
a pointed stable curve. Note the conditions
of (3.5) apply -- by assumption (3.6.2) if $r \leq 1$ and
otherwise by \cite{GKM00,0.9}. Thus by (3.5), (3.4)
holds for all $E$. 

(As foretold) we define two pointed curves $F,E$, corresponding
to closed boundary points of $\mgn$ to be 
$R$-equivalent if
the polarised varieties $Z_{P_F}$ and $Z_{P_E}$
are isomorphic, and $q_F([F]) = q_E([E])$
in $Q_E =Q_F$
under the identification of (3.1), and the map 
of (3.3) (note this makes sense since
$L_{P_X}$ is semi-ample for any 
proper closed stratum). We define a closed interior point to be
$R$-equivalent to itself, but not to any other
point. $R$ is clearly an equivalence relation
on the set of closed points.

Choose a product decomposition for each closed stratum and
let $\cP$ be their disjoint union. Let 
$\tcX$ be the disjoint union of the normalisations of
all closed strata, and $\cX^0 \subset \tcX$
the disjoint union of all (open)
strata.
Let $\cQ$ be the
disjoint union of $\mgn$ together with
the disjoint union of $Q_X$, one for each isomorphism
class of polarised varieties $Z_{P_X}$, for
$X$ a proper closed stratum. The maps $q_{P_X}$ and
$q_X$ of (3.3) give maps 
$$
\align
q_{\cP} &: \cP \rightarrow \cQ \\
q_{\tcX} &: \tcX \rightarrow \cQ \\
q        &: \cX^o \rightarrow \cQ
\endalign
$$
with $q_{\cP}$ the composition of
$q_{\tcX}$ with $\cP \rightarrow \tcX$.
$L|_{\cP}$ is pulled back from $\cQ$.

Let $\oR$ be the image of
$$
\cP \times_{\cQ} \cP \rightarrow \mgn \times \mgn. 
$$
Since $q_{\cP}$ factors through $q_{\tcX}$, 
$\oR$ is also the image of
$$
\tcX \times_{\cQ} \tcX \rightarrow \mgn \times \mgn.
$$
$R$ is by construction the set of closed points of
$$
\cX^o\times_{\cQ} \cX^o \subset \mgn \times \mgn
$$
and so contained in the closed points of $\oR$. 
$\oR$ is projective
so $\oR \subset \mgn \times \mgn$ is closed
(we will not need to know whether or not it's 
an equivalence relation).
Since $L|_{\cP}$ is pulled back from $\cQ$ and
$\cP \rightarrow \mgn$ is finite, it follows that
$p_2^*(L)$ is numerically trivial on fibres of the
first
projection $\oR \rightarrow \mgn$: Any fibre of
the first projection $\oR \rightarrow \mgn$ is the
union of images
of finitely many fibres of the first projection
$\cP \times_{\cQ} \cP \rightarrow \cP$, and thus
of finitely many fibres of $\cP \rightarrow Q$. $L_{\cP}$
is numerically trivial on fibres of $\cP \rightarrow Q$,
since it is pulled back from $Q$. 

Thus by (2.4) 
it's enough to show $L$ equivalence is
a refinement of $R$  (for then $\oR$ is
related to $L$) and for this, since $R$ is 
an equivalence relation, it's enough to 
show that for  
$Y \subset \mgn$ an irreducible curve, with 
$L \cdot Y = 0$, any two points of $Y$ are $R$-equivalent.
$Y \subset \ex L. \subset \partial \mgn.$
$Y$ meets the interior of a unique closed
stratum. Clearly it's enough to consider a point
$[E] \in Y$ in the interior of this stratum and an 
arbitrary point $[F] \in Y$. $F$ is a degeneration
of $E$ as in (3.0.1).
Let $Y' \subset P_E$
be an irreducible curve lifting $Y$, and 
$[\tF],[\tE] \in Y'$ points mapping to $[F],[E]$.
Let $p: P_E \rightarrow W$ be the map of
(3.4). $L_{P_E} \cdot Y' =0$, and $Y'$ meets
the interior of $P_E$, so by (3.4.2), $p(Y')$ is
a point, and in particular 
$p([\tF])=p([\tE])$. So by (3.4.3), $F$
and $E$ are $R$-equivalent. \qed \enddemo

%
%

\remark{3.7 Remark} If one has (3.6.2)
for all $r$, then (3.6) holds also for $p=2$ --
in the proof the assumption $p > 2$ is used only
to obtain (3.6.2) (via \cite{GKM00,0.9}) for
$r > 1$. \endremark

\demo{Proof of (0.5)} We use (3.7). So we need
to check that the
conditions of (3.6.2) hold for each of the line
bundles in question, and all $r \geq 0$. For $\psi_i$ apply 
\cite{Keel99,4.9} (where 
$\psi_i$ on $\vmgn g. {n+1}.$ is
denoted $L_{g,n}$). For $D_{g,n}$ apply
\cite{GKM00,4.7-8}. For $\lambda$: $\lambda$ is
trivial on $\vmgn 0. n.$. $\lambda$ on $\vmgn g.n.$ 
for $g \geq 2$ (resp. for $g=1$) is the
pullback of $\lambda$ from $\mg$ (resp. 
$\vmgn 1.1.$). Finally 
$\ex {\lambda_{\mg}}. \subset \partial \mg$, which
follows e.g. from Mumford's formula 
$12 \lambda = \kappa + \delta$
and the well known ampleness of $\kappa$. \qed \enddemo

\subhead \S 4 Semi-ample and EWM \endsubhead
{\bf 4.0 Assumptions:} Throughout this section $k$ is 
the algebraic
closure of a finite field. $L$ is a nef line bundle. 
All spaces are proper algebraic spaces over $k$.
Throughout we work with $\bQ$-line bundles. Note by
the base field assumption, that a $\bQ$-line
bundle is the same as its Chern class in the Neron-Severi
group.

\proclaim{4.1 Theorem} Let $L$ be a nef line bundle on 
an algebraic space $X$ proper over $k$, 
the algebraic closure of
a finite field. Then $L$ is semi-ample iff $L$ is EWM. 
\endproclaim

\remark{4.1.1 Remark} Note (4.1) is equivalent
to the statement that the class $c_1(L)$ in
the Neron-Severi group is pulled back
from $Z$ for $g:X \rightarrow Z$ the associated map.
This formulation looks at first reasonable for any
field, but unfortunately fails: 
Take for example a smooth curve $C \subset X$
of negative self intersection on a projective surface
so that the contraction, $g$, of $C$ is not projective.
Take $D$ an effective combination of $C$ and ample $H$
so that $D \cdot C = 0$. Then $g$ is the map associated
to $D$, 
but 
the first Chern class (in the Neron-Severi group)
 cannot descend, as it would
represent an ample divisor on the image. I suspect
there are also counter-examples with $Z$ projective. \endremark

We use the following in proving (4.1):

\proclaim{4.2 Lemma} Let $g: X \rightarrow Z$ be 
a proper surjection with geometrically connected fibres,
with $X$ projective.
Assume $L$ is nef, $EWM$ and $g$-numerically trivial. 
Assume there is a proper $Y \rightarrow X$, with
$Y \rightarrow Z$ surjective and a 
line bundle $L_Z$ so that $L|_Y = L_Z|_Y$.
Then $L = L_Z|_X$ (as $\bQ$-line bundles). \endproclaim
\demo{Proof} 
It's enough by (1.7) to show
$L$ is pulled back from $Z$. As we work 
in the Neron-Severi group, to check
$L = {L_Z}|_X$ it's enough to check the
intersection number of either side 
with a given subcurve $C \subset X$.
We can replace $Z$ by the image of $C$ and everything else
by the obvious pullbacks. 
We can replace $X$ by a general hyperplane 
section containing $C$, 
so long as $g$ has
fibre dimension at least $2$ (for then the restriction of
$g$ to the hyperplane section still has connected fibres). 

So we can assume $X$ has dimension at
most two  and $Z$ is a curve.
It's enough, by (1.7), to show $L$ is semi-ample, for then,
since $g$ has connected fibres, the associated map 
will the composition of $g$ with a finite universal
homeomorphism. 
Since $L$ is $EWM$, $L$-equivalence is
bounded, so 
by (0.1) to check semi-ampleness we can replace $X$ by a 
connected component of a desingularisation.
Any nef line bundle on a curve over a
finite field is semi-ample 
so the next result applies:
\qed \enddemo

For the next lemma we allow arbitrary base field.

\proclaim{4.4 Lemma} Let $p:S \rightarrow C$ be 
a map from an irreducible non-singular projective surface to 
a normal curve. Let $L$ be a nef, $p$-numerically
trivial line bundle. The class $c_1(L)$ in
the Neron-Severi group is pulled back from $C$.
\endproclaim
\demo{Proof} We can assume $p$ is surjective, otherwise
$L$ is numerically trivial and the result is obvious. 
Then apply the next lemma with $\beta$ the class
of a general fibre and $\gamma = c_1(L)$. 
\qed \enddemo

\proclaim{4.5 Lemma} Let $S$ be a smooth
projective surface and $\beta \in \nsd S.$ such that
$\beta^2 = 0$ and $\beta \cdot H > 0$, for some
ample class $H$. If $\gamma \cdot \beta = 0$
then $\gamma^2 \leq 0$, with equality iff
$\gamma$ is equivalent to a multiple of $\beta$. 
\endproclaim
\demo{Proof} This is immediate from the hodge
index theorem. Indeed suppose $\gamma \cdot \beta =0$.
Then 
$$
\align
(\gamma - \frac{\gamma \cdot H}{\beta \cdot H} \beta) \cdot H &= 0 \\
(\gamma - \frac{\gamma \cdot H}{\beta \cdot H} \beta)^2 & = \gamma^2. \qed
\endalign
$$
\enddemo
\demo{Proof of 4.1} We induct on
the dimension of $X$. Since $L$ is EWM, $L$-equivalence
is algebraic. Thus by (0.1) we may replace $X$ by
anything of the same dimension that surjects properly
onto it, and so may assume $X$ is normal, irreducible, and 
projective. Let $g: X \rightarrow Z$ be
the map associated to $L$. 

Suppose first the relative dimension is at least $2$.
Let $Y \subset X$ be an ample hypersurface. 
$Y \rightarrow Z$ is a surjection with geometrically
connected
fibres, the Stein factorisation $Y \rightarrow Z_Y$
is the map associated to $L_Y$ and $Z_Y \rightarrow Z$
is a finite universal homeomorphism. 
Thus by induction and \cite{Keel99,1.4}, 
$L_Y$ is pulled back from $Z$, necessarily
from an ample line bundle. Now apply (4.2). 

Suppose next that $g$ is birational. Then $L$ is
big. $L|_{\ex L.}$ is semi-ample by induction, so $L$
is semi-ample by \cite{Keel99,0.2}. 

So we can assume  $g$ has relative dimension one. 

Let $Y' \subset X$
be a general hyperplane section, and let $Y \rightarrow Y'$
be a finite cover so the (finite) field extension
$K(Z) \subset K(Y)$ is normal, with group 
Galois group $G$. 
Let $g_Y: Y \rightarrow Z_Y$ be the Stein factorisation
of $Y \rightarrow Z$. Note it is birational,
and precisely 
the map associated
to $L_Y$.  By induction, $L_Y$ is pulled back from
ample $L_{Z_Y}$. $G$ acts
on $(Z_Y,L_{Z_Y})$. Since 
$g_Y$ is birational, $G$ is the Galois group
of the normal field extension $K(Z_Y)/K(Z)$.
Thus, since $Z_Y \rightarrow Z$ is finite,
the quotient $Z_Y/G$ is a finite purely
inseparable cover of $Z$. 
By (4.2) it's enough to show the $\bQ$-line
bundle 
$L_{Z_Y}$ is pulled back from $Z$ and so,
by \cite{Keel99,1.4},
it's enough to show
$c_1(L_Y)$ is preserved by $G$. 
Thus we need only check 
$$
\gamma^*(L_{Z_Y}) \cdot C = L_{Z_Y} \cdot C \tag{4.6}
$$
for a given $\gamma \in G$ and irreducible $C \subset Z_Y$. 
Let $f$ be the map $Z_Y \rightarrow Z$ and $D = f(C)$. 
Note for (4.6) it's enough (since $f$ is $G$ invariant),
to show $L_{Z_Y}|_{f^{-1}(D)}$ is pulled
back from $D$ and for this, by (1.7), it's enough to
show $L|_{g^{-1}(D)}$ is pulled back from $D$. Since
$g$ (and so $g|_{g^{-1}(D)}$) has geometrically
connected fibres, it's enough, by \cite{Keel99,1.4}, to
show $L|_{g^{-1}(D)}$ is semi-ample (for then
the associated map will be the composition of
(the restriction of) $g$ with a finite universal
homeomorphism). Thus we can replace $X$ by
$g^{-1}(D)$, and so reduce to the case when $Z$ (the
image of the associated map) is a curve. 
Now return to the beginning of
the proof and retrace the argument. We reduce
to the case $X$ an irreducible normal surface
and $Z$ a normal curve. We can replace $X$ by a desingularisation, 
and apply (4.4). \qed \enddemo

\subhead \S 5 counter example \endsubhead

We begin with a general construction. Let
$f_i: X \rightarrow X_i$, $L_{X_i}$ be as in (0.12)
with $L_{X_i}$ ample.  Let $Y$ be obtained
by gluing $Y_1 := X \times X_1$ to $Y_2:= X \times X_2$
along the graphs of $f_1,f_2$. The line bundles
$L_{Y_i} : = \pi_{X_i}^*(L_{X_i})$ are obviously 
semi-ample and glue to give a nef line bundle
$L_Y$ on $Y$, with ${L_Y}|_X = L_X$. 

\proclaim{5.1 Lemma} Notation as above. Assume $X$ is
geometrically connected. 
If $L_Y$
is EWM (resp semi-ample), then a finite pushout (resp.
polarised pushout) exists
for $f_1,f_2$ (resp $(f_i,L_{X_i})$).  \endproclaim
\demo{Proof} 
Let $g: Y \rightarrow Z$ be the map associated to
$L_Y$. Note $\pi_{X_i}$ is the map associated 
to $L_{Y_i}$ thus $g|_{Y_i}$ factors through 
$\pi_{X_i}$. In particular the finite map
$g|_X: X \rightarrow Z$ factors through $f_1$ 
and through $f_2$, and is thus a pushout for $f_1,f_2$.
If
$L_Y$ is pulled back from $L_Z$, then the pushout
is polarised. \qed \enddemo

\proclaim{5.2 Lemma (Koll\'ar)} Let $k$ be the algebraic
closure of $\bF_p$. The intersection of 
the subfields $k(T^p + T^{p-1})$ and $k(T^{p-1})$ of the
polynomial field $k(T)$ is $k$. In particular the
corresponding pair of endomorphisms of $\bP^1$ has
no finite pushout. \endproclaim
\demo{Proof} Suppose on the contrary
the fields have nontrivial intersection. Then
there exist monic polynomials $f_i,g_i$,
with 
$$
(f_1(T),f_2(T)) = (g_1(T),g_2(T)) = 1.
$$
such that 
$$
f_1(T^p + T^{p-1})/f_2(T^p + T^{p-1}) = 
g_1(T^{p-1})/g_2(T^{p-1}) \tag{5.2.1}
$$
with $f_1,g_1$ non-constant. 
There exist $a_1,a_2$ with 
$$
a_1(T) f_1(T) + a_2(T) f_2(T) =1.
$$
Replacing $T$ by $T^p + T^{p-1}$ we conclude
the $f_i(T^p + T^{p-1})$ are relatively prime.
Similarly the $g_i(T^{p-1})$ are relatively prime.
Now (5.2.1) implies 
$f_i(T^p + T^{p-1}) = g_i(T^{p-1})$. Thus we
can assume $f_2=g_2 =1$. 
Now choose $f$ and $g$ of
minimal positive degree so that 
$f(T^p + T^{p-1}) - g(T^{p-1})$ is of
degree at most one, say $aT + b$. We
can absorb $b$ into $f$, and assume $b =0$.
Differentiating (and using characteristic $p$)
we see that 
$$(p-1)T^{p-2}(f'(T^p + T^{p-1}) -g'(T^p)) = a
$$
is constant.  We conclude that 
$a =0$, and so (by minimality of degree)
$f' = g' = c$ is constant. Thus for
some polynomials $F$ and $G$,
$f(T) = F(T)^p + c T$, $g(T) = G(T)^p + c T$, and we have
$$
F(T^p + T^{p-1})^p + c T^p = G(T^{p-1})^p
$$
But then 
$$
F(T^p + T^{p-1}) - G(T^{p-1})
$$
is linear, contradicting the minimality. \qed \enddemo

Now take $X=X_1=X_2 = \bP^1$ and $L_{X_i} = \ring .(1)$,
and $f_1,f_2$ the endomorphisms of (5.2). 
The construction of (5.1) gives a 
non semi-ample 
line bundle on a surface defined over a finite field 
whose restriction to the normalisation is semi-ample.

\end